\documentclass{article}
\usepackage{graphicx} 
\usepackage[a4paper, margin=3.14cm]{geometry}
\usepackage{amsmath}
\usepackage{mathtools}
\usepackage{amsfonts}
\usepackage{amssymb}
\usepackage{amsthm}
\usepackage{tikz} 
\usetikzlibrary{decorations.pathreplacing}
\usepackage[vcentermath]{youngtab}
\usetikzlibrary{3d}
\usetikzlibrary{perspective}
\usepackage{float}
\usepackage{comment}
\usepackage{setspace}

\theoremstyle{plain}
\newtheorem{theorem}{\textbf{Theorem}}[section]

\newtheorem{prop}{\textbf{Proposition}}[section]
\newtheorem{cor}{\textbf{Corollary}}[section]

\theoremstyle{definition}
\newtheorem{defn}{\textbf{Definition}}[section]

\newtheorem{example}{Example}[section]

\theoremstyle{remark}
\newtheorem{remark}{Remark}[section]

\DeclarePairedDelimiter{\floor}{\lfloor}{\rfloor}
\DeclarePairedDelimiter{\brac}{ \{}{ \}}
\DeclarePairedDelimiter{\hk}{  (}{ )}
\usetikzlibrary{math}

\newcommand{\modgram}[3]{%
  \begin{tikzpicture}[scale=0.5]
    \foreach \part [count=\i from 0] in {#2} {
      
      \pgfmathtruncatemacro{\q}{\part / #1}
      \pgfmathtruncatemacro{\r}{mod(\part, #1)}
      
      \ifnum\q>0
        \foreach \col in {1,...,\q} {
          \pgfmathtruncatemacro{\x}{\col - 1}
          \pgfmathtruncatemacro{\y}{-\i}
          
          \draw (\x, \y) rectangle (\x + 1, \y + 1);
          \node at (\x + 0.5, \y + 0.5) {$#1$};
        }
      \fi
      
      \ifnum\r>0
        \pgfmathtruncatemacro{\x}{\q}
        \pgfmathtruncatemacro{\y}{-\i}
        
        \draw (\x, \y) rectangle (\x + 1, \y + 1);
        \node at (\x + 0.5, \y + 0.5) {$\r$};
      \fi
    }
    
    \ifnum#3>0
      \fill[blue!20, opacity=0.4] (0, -#3+1) rectangle (#3, 0+1);
      \draw[blue, thick] (0, -#3+1) rectangle (#3, 0+1);
    \fi
  \end{tikzpicture}%
}

\title{Durfee Cubes and Higher-Order $q$-Identities for Partitions}
\author{Marthinus Michael Dreeckmeier\\
Darlison Nyirenda\\
\small School of Mathematics\\
\small University of the Witwatersrand\\
\small P.~O. Wits 2050, Johannesburg, South Africa
}
\date{\today}

\newcommand{\cube}[3]{
  \begin{scope}[shift={(#1,#2,#3)}]
    \draw[fill=gray!5, opacity=0.6] (0,0,0) -- (1,0,0) -- (1,0,1) -- (0,0,1) -- cycle; 
    \draw[fill=gray!5, opacity=0.6] (0,0,0) -- (0,1,0) -- (0,1,1) -- (0,0,1) -- cycle; 
     \draw[fill=gray!5, opacity=0.6] (0,0,0) -- (0,1,0) -- (1,1,0) -- (1,0,0) -- cycle; 

    \draw[fill=gray!50, opacity=0.89] (0,0,1) -- (1,0,1) -- (1,1,1) -- (0,1,1) -- cycle; 
    \draw[fill=gray!30, opacity=0.89] (1,0,0) -- (1,1,0) -- (1,1,1) -- (1,0,1) -- cycle; 
    \draw[fill=gray!15, opacity=0.89] (0,1,0) -- (1,1,0) -- (1,1,1) -- (0,1,1) -- cycle; 
  \end{scope}
}

\begin{document}

\maketitle
\begin{abstract}
     In the study of $q$-series, there do not seem to be any product-sum $q$-identities with an exponent in the summation of order greater than 2. We generalise the notion of the Durfee square for partitions by introducing higher-dimensional Durfee cubes. We use this new concept to derive higher-order generalisations of classical product-sum $q$-identities for certain classes of restricted and unrestricted partitions.
\end{abstract}

\section{Introduction}
A \textit{partition} of an integer $n$, given by $\lambda=(\lambda_1,\lambda_2,\dots,\lambda_s)$, is a one-dimensional array of positive integers that satisfies $\lambda_i\geq \lambda_{i+1}$ for $1\leq i < s$ and
$$
\sum_{1\leq i\leq s}\lambda_i=n.
$$
A product-sum identity is any identity that relates a summation to a product form. A simple example of this is the geometric series formula. For $|q|<1$, we have
$$
\sum_{s\geq 0}q^s=\frac{1}{1-q}.
$$
Many product-sum identities for partitions can be derived via Durfee squares \cite{r4}. For unrestricted partitions, we have the Cauchy identity. Before we give the identity, we will first introduce some notation. The \textit{$q$-Pochhammer symbol} $(a;q)_s$ is defined as follows
$$
(a;q)_s:=\prod_{i=0}^{s-1}(1-aq^i).
$$
When $s=0$, we have an empty product, hence $(a;q)_0=1$.
\begin{theorem}[Cauchy]\label{thmp1}
    For $|q|<1$, we have
\begin{equation}\label{eqintro1}
    \frac{1}{(q;q)_\infty}=\sum_{s\geq 0}\frac{q^{s^2}}{(q;q)_s^2}.
\end{equation}
\end{theorem}
For partitions into distinct parts and partitions into distinct odd parts, similar product-sum identities can be derived with the help of Durfee squares, where the exponent of $q$ in the summation is quadratic. For partitions into distinct parts, we have
\begin{equation}\label{ieq2}
    (-zq;q)_\infty=\sum_{n\geq0}z^nq^{n(3n+1)/2}(zq^{2n+1}+1)\frac{(-zq;q)_n}{(q;q)_n}.
\end{equation}
Notice that Euler's pentagonal number theorem can be derived from equation (\ref{ieq2}).
For partitions into distinct odd parts, we have
\begin{equation}\label{ieq3}
    (-q;q^2)_\infty=\sum_{n\geq0}\frac{q^{n^2}}{(q^2;q^2)_n}.
\end{equation}

M.~D.~Hirschhorn \cite{r1} gave a proof of Jacobi's four-square theorem and Jacobi's two-square theorem \cite{twosquare} using only the Jacobi triple product identity and basic $q$-series manipulations. The key component of the proof was the fact that the Jacobi triple product identity has a square in the exponent of $q$ in the summation. However, attempting to extend this strategy for sums of cubes or sums of higher powers does not seem possible. This is because it seems that there is no cubic or higher-order generalisation of the Jacobi triple product identity. In this study, we will not be solving problems regarding sums of cubes, nor will we find a generalisation of the Jacobi triple product. Instead, we will address what we believe to be a gap in the literature highlighted by this: there do not seem to be, as far as we are aware, any known product-sum $q$-identities that have exponents of higher order than quadratic in the summation. In particular, this would be any identity roughly of the form $$\text{product}=\sum_{s} f(s,q)q^{g(s)},$$ where $g(s)$ is a cubic or higher-order expression and $f(s,q)$ could consist of products and other terms with respect to $q$ and $s$.

Generalisations of Durfee squares have led to several extensions of classical partition identities \cite{r2,r3}. These developments suggest that modifying the Durfee decomposition can reveal new $q$-series structures. Motivated by this perspective and our previously stated motivation, we set out to generalise the Durfee square to Durfee cubes and Durfee hypercubes. For example, we show in Section \ref{section2} that, with the aid of Durfee hypercubes, equation (\ref{eqintro1}) can be generalised to the following identity:
\begin{align}
    \frac{1}{(zq;q)_{\infty}}=&1+\sum_{s \geq 1}\sum_{r\geq0}\frac{z^sq^{s^k+rs}}{(q;q)_{s-1}(zq;q)_{s^{k-1}+\alpha(r,s)}}\notag\\
    =&1+\sum_{s \geq 1}\sum_{r\geq0}\frac{z^{s^{k-1}+r}q^{s^k+rs}}{(zq;q)_{s-1}(q;q)_{s^{k-1}+\alpha(r,s)}},\notag\\
\end{align}
where $\alpha(r,s)=\min{\{r,(s+1)^{k-1}-s^{k-1}-1\}}$.
Notice that setting $k=2$ in the above identity results in the Cauchy identity.
We also generalise the identities in equations (\ref{ieq2}) and (\ref{ieq3}).
\section{Preliminaries, Definitions, and Notation}
In this section, we state known definitions and theorems that we use or generalise in our results. Most of these can be found in the book \textit{The theory of partitions} by George Andrews \cite{r4}, or the book \textit{Integer partitions} \cite{ip}.

We will use the symbol $\mathbb{N}_0$ to denote the set of natural numbers including 0. Throughout this work, $n$ shall stand for a nonnegative integer.

Partitions can be represented visually via Ferrers diagrams \cite{ip}. To obtain the \textit{Ferrers diagram} of a partition $\lambda=(\lambda_1,\lambda_2,\dots,\lambda_s)$, we place dots or squares in rows and columns where the $i$th row of the diagram contains $\lambda_i$ dots. Another visual representation of a partition is called the $m$-modular Ferrers diagram of a partition \cite{macmahon1923}. To obtain the \textit{$m$-modular Ferrers diagram} of a partition $\lambda=(\lambda_1,\lambda_2,\dots,\lambda_s)$, write $\lambda_i=mq_i+r_i$, where $1<r_i\leq m$, then the $i$th row in the diagram will contain $q_i$ squares/nodes with weight $m$ and one square/node with weight $r_i$. The \textit{Durfee square} of a partition is the largest square that can fit in the Ferrers diagram of that partition.
For our purposes, we give a more formal definition for the Durfee square.
\begin{defn}\label{defp1}
    A partition of $n$ contains a \textit{square} of side length $s$ if the partition contains at least $s$ parts greater than or equal to $s$.
\end{defn}
\begin{defn}\label{defp2}
    The \textit{Durfee square} of a partition of $n$ has a side length of $s$ if the partition contains a square of side length $s$ but not a square of side length $s+1$.
\end{defn}
The following proposition can be used to obtain the side length of the Durfee square of a partition.
\begin{prop}\label{lmap1}
        Let $\lambda=(\lambda_1,\,\lambda_2,\dots,\lambda_j)$ be a partition of $n$ with $\lambda_1\geq\lambda_1\geq\cdots\geq\lambda_j$. Then the side length, $s$, of the Durfee square can be determined by
    \begin{equation*}
        s=\vert \{\lambda_i:1\leq i\leq j \text{ and } \lambda_i\geq i\}\vert.
    \end{equation*}
\end{prop}
We generalise the following theorems in Section \ref{section2}.
\begin{theorem}\label{thmp3}
    For $|q|<1$, we have
    $$
    \frac{1}{(q;q)_\infty}=\sum_{s\geq 0}\frac{q^{s}}{(q;q)_s}.
    $$
\end{theorem}
\begin{theorem}\label{thmc3.4}
    For $|q|<1$, we have
    $$\sum_{j=0}^sq^{j^2}{\begin{bmatrix}
        s\\
        j
    \end{bmatrix}}_q^2={\begin{bmatrix}
        2s\\
        s
    \end{bmatrix}}_q.$$
\end{theorem}
Let $\lambda$ be a partition of $n$. The \textit{conjugate} of $\lambda$, denoted by $\lambda'$, is the partition obtained by reading the columns of the Ferrers diagram of $\lambda$ as rows and rows as columns.

\section{Identities Arising from Generalised Durfee Hypercubes}\label{section2}
Most known $q$-series identities contain terms such as $q^n$ or $q^{n^2}$ in the summations. In this section, we will generalise Durfee squares and use them to find $q$-identities with terms of the form $q^{n^k}$ in the summation for values of $k>2$.
\subsection{Introductory Constructions and Definitions}
\begin{figure}[H]
    \centering
    \begin{tikzpicture}[scale=0.7]

\coordinate (0) at (0,0);

\draw[->, thick] (0) -- (3,0) node[right] {$x$};
\draw[->, thick] (0) -- (-2*0.85,-1.5*0.85) node[below left] {$y$};
\draw[->, thick] (0) -- (0,3) node[above] {$z$};

\end{tikzpicture}
    \caption{The $xyz$-coordinate system.}
    \label{fig3.1}
\end{figure}
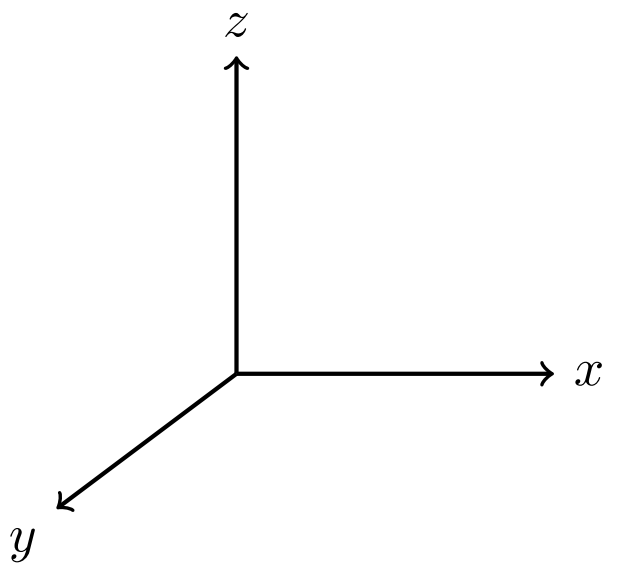

For three-dimensional visualisations, we use an $xyz$-coordinate system where the $x$-axis stretches from left to right, the $z$-axis points upwards, and the $y$-axis points out of the page (see Figure \ref{fig3.1}).
We introduce the concept of stacked Ferrers diagrams.
\begin{defn}
    An \textit{$m$-stacked Ferrers diagram} of a partition $\lambda=(\lambda_1,\,\lambda_2,\dots,\lambda_j)$ of $n$ is obtained by stacking the first part $\lambda_1$ starting at $(0,0,0)$ extending in the positive $x$ direction. The next part $\lambda_2$ is stacked on top of $\lambda_1$. Then $\lambda_3$ is stacked on $\lambda_4,\dots$, finally $\lambda_m$ is stacked on $\lambda_{m-1}$. Once we reach this point, stack the next part, $\lambda_{m+1}$, at the bottom next to $\lambda_1$ (with the shift happening in the $y$-axis). Then $\lambda_{m+2}$ is stacked on $\lambda_{m+1}$ and so on up to $\lambda_{2m}$ on $\lambda_{2m-1}$. Repeat this process until all parts are depleted, and go up at most $m$ parts before shifting along the $y$-axis.
\end{defn}
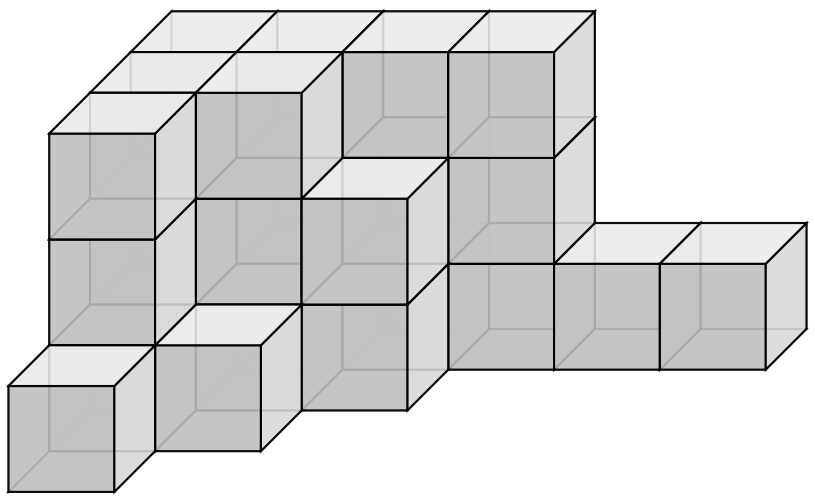
\begin{figure}[H]
    \centering
    \begin{tikzpicture}[rotate around x=0, rotate around y=0, rotate around z=0, scale=0.7]
    \cube{0}{0}{0} \cube{1}{0}{0} \cube{2}{0}{0} \cube{3}{0}{0} \cube{4}{0}{0} \cube{5}{0}{0}
    
    \cube{0}{1}{0} \cube{1}{1}{0} \cube{2}{1}{0} \cube{3}{1}{0}

    \cube{0}{2}{0} \cube{1}{2}{0} \cube{2}{2}{0} \cube{3}{2}{0}

    \cube{0}{0}{1} \cube{1}{0}{1} \cube{2}{0}{1}

    \cube{0}{1}{1} \cube{1}{1}{1} \cube{2}{1}{1}

    \cube{0}{2}{1} \cube{1}{2}{1}

    \cube{0}{0}{2} \cube{1}{0}{2}

    \cube{0}{1}{2}

    \cube{0}{2}{2}

    \cube{0}{0}{3}

\end{tikzpicture}
    \caption{The 3-stacked Ferrers Diagram for $(6,4,4,3,3,2,2,1,1,1)$.}
\end{figure}
Note that the three-dimensional visualisation of stacked Ferrers diagrams can be collapsed into two dimensions. This can be done by keeping the top view shape in the $xy$-plane and writing the height of the blocks in the $z$ direction inside each square. As an example, consider the 2-stacked diagrams for the partition $\tau= (5,\,4,\,3,\,3,\,3,\,2,\,1)$ of 21.

\begin{figure}[H]
    \centering
    \begin{tikzpicture}[rotate around x=0, rotate around y=0, rotate around z=0, scale=0.7]
    \cube{0}{0}{0} \cube{1}{0}{0} \cube{2}{0}{0} \cube{3}{0}{0} \cube{4}{0}{0}
    
    \cube{0}{1}{0} \cube{1}{1}{0} \cube{2}{1}{0} \cube{3}{1}{0} 

    \cube{0}{0}{1} \cube{1}{0}{1} \cube{2}{0}{1}

    \cube{0}{1}{1} \cube{1}{1}{1} \cube{2}{1}{1}

    \cube{0}{0}{2} \cube{1}{0}{2} \cube{2}{0}{2}

    \cube{0}{1}{2} \cube{1}{1}{2}

    \cube{0}{0}{3}

\end{tikzpicture}%
\qquad
\begin{tikzpicture}[scale=0.7]

\draw (0,0) rectangle (1,1);
\node at (0.5,0.5) {$2$};

\draw (1,0) rectangle (2,1);
\node at (1.5,0.5) {$2$};

\draw (2,0) rectangle (3,1);
\node at (2.5,0.5) {$2$};

\draw (3,0) rectangle (4,1);
\node at (3.5,0.5) {$2$};

\draw (4,0) rectangle (5,1);
\node at (4.5,0.5) {$1$};

\draw (0,-1) rectangle (1,0);
\node at (0.5,-0.5) {$2$};

\draw (1,-1) rectangle (2,0);
\node at (1.5,-0.5) {$2$};

\draw (2,-1) rectangle (3,0);
\node at (2.5,-0.5) {$2$};

\draw (0,-2) rectangle (1,-1);
\node at (0.5,-1.5) {$2$};

\draw (1,-2) rectangle (2,-1);
\node at (1.5,-1.5) {$2$};

\draw (2,-2) rectangle (3,-1);
\node at (2.5,-1.5) {$1$};

\draw (0,-3) rectangle (1,-2);
\node at (0.5,-2.5) {$1$};

\end{tikzpicture}
 \caption{The 2-stacked Ferrers diagram for $\tau$ and its 2-dimensional visualisation.}
 \label{fig:1}
\end{figure}
Also note that modular Ferrers diagrams can be viewed in three dimensions, where the number in each square represents the number of blocks stacked in the positive $z$ direction.
\begin{remark}\label{remr1}
    The $m$-stacked and $m$-modular Ferrers diagrams are conjugates of each other.
\end{remark}
\begin{example}
    The conjugate of the partition $\tau= (5,\,4,\,3,\,3,\,3,\,2,\,1)$ is $\tau'=(7,\,6,\,5,\,2,\,1)$. As can be seen in Figure \ref{fig3.3}, the 2-modular Ferrers diagram for $\tau'$ is the conjugate of the 2-stacked Ferrers diagram for $\tau$ in Figure \ref{fig:1}.
    \begin{figure}[H]
    \centering
    \begin{tikzpicture}[scale=0.5]

        \draw (0,0) rectangle (1,1); \node at (0.5,0.5) {$2$};
        \draw (1,0) rectangle (2,1); \node at (1.5,0.5) {$2$};
        \draw (2,0) rectangle (3,1); \node at (2.5,0.5) {$2$};
        \draw (3,0) rectangle (4,1); \node at (3.5,0.5) {$1$};

        \draw (0,-1) rectangle (1,0); \node at (0.5,-0.5) {$2$};
        \draw (1,-1) rectangle (2,0); \node at (1.5,-0.5) {$2$};
        \draw (2,-1) rectangle (3,0); \node at (2.5,-0.5) {$2$};

        \draw (0,-2) rectangle (1,-1); \node at (0.5,-1.5) {$2$};
        \draw (1,-2) rectangle (2,-1); \node at (1.5,-1.5) {$2$};
        \draw (2,-2) rectangle (3,-1); \node at (2.5,-1.5) {$1$};

        \draw (0,-3) rectangle (1,-2); \node at (0.5,-2.5) {$2$};

        \draw (0,-4) rectangle (1,-3); \node at (0.5,-3.5) {$1$};

    \end{tikzpicture}
    \caption{The 2-modular Ferrers diagram for $\tau'$.}
    \label{fig3.3}
\end{figure}    
\end{example}
We call a $k$-dimensional cube a \textit{$k$-cube}.
We now generalise Definitions
\ref{defp1} and \ref{defp2} for stacked and modular Ferrers diagrams.
\begin{defn}
    A partition of $n$ contains a \textit{stacked $k$-cube} of side length $s$ if the partition contains at least $s^{k-1}$ parts greater than or equal to $s$.
    A partition of $n$ contains a \textit{modular $k$-cube} of side length $s$ if the partition contains at least $s$ parts greater than or equal to $s^{k-1}$.
\end{defn}
\begin{defn}
    The \textit{Durfee stacked $k$-cube} of a partition of $n$ has a side length $s$ if the partition contains a stacked $k$-cube of side length $s$, but not a stacked $k$-cube of side length $s+1$.
    The \textit{Durfee modular $k$-cube} of a partition of $n$ has a modular $k$-cube of side length $s$, but not a modular $k$-cube of side length $s+1$.
\end{defn}
\begin{remark}
    Let $k\geq2$.
    \begin{itemize}
    \item For every partition of $n$, the Durfee modular $k$-cube is the largest modular $k$-cube and the Durfee stacked $k$-cube is the largest stacked $k$-cube.
    \item Every partition of $n$ contains a unique Durfee modular $k$-cube and a unique Durfee stacked $k$-cube.
    \end{itemize}
\end{remark}
As a consequence of Remark \ref{remr1}, we have the following corollary.
\begin{cor}
    Let $k\geq2$. The number of partitions of $n$ that contain a modular $k$-cube of side length $s$ is equal to the number of partitions of $n$ that contain a stacked $k$-cube of side length $s$.
\end{cor}
\begin{proof}
    If a partition of $n$ contains a modular $k$-cube of side length $s$, then by Remark \ref{remr1}, its conjugate contains a stacked $k$-cube of side length $s$. Because the conjugation is an involution, we have a bijection.
\end{proof}
\begin{cor}
    Let $k\geq2$. Then, the number of partitions of $n$ where the Durfee modular $k$-cube has a side length of $s$ is equal to the number of partitions of $n$ where the Durfee stacked $k$-cube has a side length of $s$.
\end{cor}
\begin{proof}
    Let $m_s(n)$ and $st_s(n)$ be the number of partitions of $n$ that contain a modular $k$-cube of side length $s$ and a stacked $k$-cube of side length $s$, respectively. By the preceding corollary, we know that $m_s(n)=st_s(n)$ and $m_{s+1}(n)=st_{s+1}(n)$. Let $M_s(n)$ and $ST_s(n)$ be the number of partitions of $n$ that contain a  Durfee modular $k$-cube of side length $s$ and a Durfee stacked $k$-cube of side length $s$, respectively. Note that $m_{s+1}(n)=m_s(n)-M_s(n)$ and $st_{s+1}(n)=st_s(n)-ST_s(n)$. Now, since $m_{s+1}(n)=st_{s+1}(n)$, we have $m_s(n)-M_s(n)=st_s(n)-ST_s(n)$. Finally, applying $m_s(n)=st_s(n)$, we obtain the desired result.
\end{proof}
We did not simply apply the conjugate again to prove the above corollary for Durfee $k$-cubes because, if you take the conjugate, it still raises the question of whether it contains a larger $k$-cube. However, keeping track of the number of parts or the length of parts shows that the conjugate does not contain a larger $k$-cube.

Next, we want to know how to determine the side length of modular and stacked Durfee cubes. Consider the following proposition, which is a generalisation of Proposition \ref{lmap1}. Note that the Durfee stacked $k$-cube and Durfee modular $k$-cube of a partition need not have the same side length.

\begin{prop}
    Let $k\geq2$. Let $\lambda=(\lambda_1,\,\lambda_2,\dots,\lambda_j)$ be a partition of $n$ in non-increasing order. Let the Durfee stacked $k$-cube of $\lambda$ have a side length of $s_1$ and the Durfee modular $k$-cube of $\lambda$ have a side length of $s_2$. Then
\begin{align}
    s_1&=\left\lfloor\left(\vert\{\lambda_i\geq i^{1/{(k-1)}}:1\leq i\leq j\}\vert\right)^{1/(k-1)}\right\rfloor,\label{eq3.1}\\
    s_2&=\vert\{\lambda_i\geq i^{{(k-1)}}:1\leq i\leq j\}\vert.\label{eq3.2}
\end{align}
\end{prop}
\begin{proof}
Let $\lambda$ and $s_1,\,s_2$ be as in the conditions of the proposition. Note that equation (\ref{eq3.2}) follows directly from the definition of the Durfee modular $k$-cube. We now show that equation (\ref{eq3.1}) holds. Clearly, once $i=s_1^{k-1}$, then you have $\lambda_i\geq s_1$ (since there are at least $s_i^{k-1}$ terms greater than or equal to $s$. Applying the same logic to $s_i+1$ to find an upper bound, we obtain
$$
(s_i)^{k-1}\leq \vert\{\lambda_i\geq i^{1/{(k-1)}}:1\leq i\leq j\}\vert < (s_i+1)^{k-1}.
$$
Taking the $(k-1)$th root
$$
s_i\leq \hk*{\vert\{\lambda_i\geq i^{1/{(k-1)}}:1\leq i\leq j\}\vert}^{1/(k-1)} < s_i+1.
$$
Finally, taking the floor gives the desired result.
\end{proof}
\begin{example}\label{ex3.2}
We use our formula to determine the side length of the Durfee stacked $3$-cube of $\tau= (5,\,4,\,3,\,3,\,3,\,2,\,1)$. We have $5\geq1^{1/2}$, $4\geq2^{1/2}$, $3\geq3^{1/2}$, $3\geq4^{1/2}$, $3\geq5^{1/2}$, but $2\not\geq 6^{1/2}$. Therefore, we have 5 terms satisfying the condition in equation (\ref{eq3.1}). Hence, $s_1=\left\lfloor5^{1/2}\right\rfloor=2$. 

We now determine the side length $s_2$ of the Durfee modular $3$-cube of $\tau'=(7,\,6,\,5,\,2,\,1)$ using equation (\ref{eq3.2}). We have $7\geq1^2$, $6\geq2^2$, but $5\not\geq 3^{1/2}$. Hence, $s_2=2$.  
\end{example}

\subsection{New Product-Sum Identities Arising from\\ Durfee Hypercubes}
Building on the diagrams and new definitions in the previous subsection, we generalise some of the classical product-sum identities discussed in the introduction. Before we jump into the general cases, we first prove identities using three-dimensional Durfee cubes to build a strong foundation and intuitive understanding of the technique.
\begin{theorem}\label{thmr1}
For $|q|<1$, we have
\begin{align}\label{eq2.1}
    \frac{1}{(q;q)_{\infty}}=1+\sum_{s \geq 1}\sum_{r\geq0}\frac{q^{s^3+rs}}{(q;q)_{s-1}(q;q)_{s^2+\alpha(r,s)}}, \,\text{where } \alpha(r,s)=\min{\{r,2s\}}.
\end{align}
\end{theorem}
\begin{proof}
Consider a partition $\lambda=(\lambda_1,\,\lambda_2,\dots,\lambda_j)$ of $n$ with a Durfee stacked $3$-cube of side length $s$. Consider the $s$-stacked Ferrers diagram for $\lambda$. We demonstrate this proof with an example. Consider the partition of 21: $\tau= (5,\,4,\,3,\,3,\,3,\,2,\,1)$. Recall from Example \ref{ex3.2} that the largest cube for $\tau$ has a side length of $2$.
Break the Ferrers diagram for $\lambda$ into four separate sections/partitions. The first is simply the cube of side length $s$, which contributes the term $q^{s^3}$. The second consists of at most $s^2$ parts to the right of the right face of the cube. The third consists of parts that are less than or equal to $s$ in front of the front face of the cube. The fourth section consists of the blocks that remain and are in front of the second and to the right of the third sections.
As an example, consider Figure \ref{fig:2}.
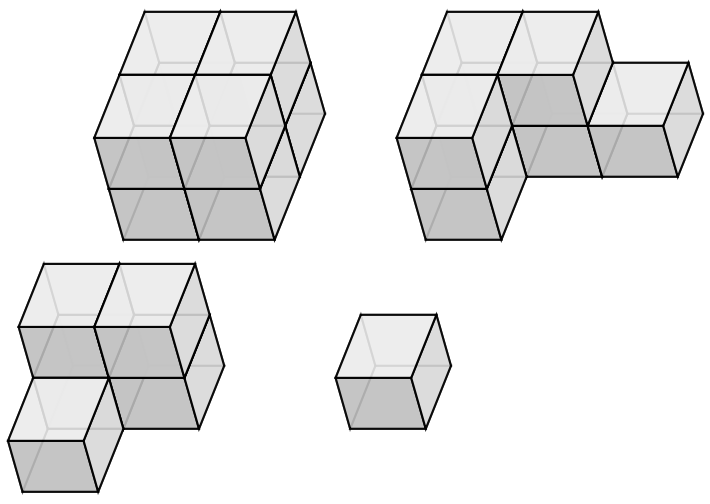
\begin{figure}[H]
    \centering
    \begin{tikzpicture}[rotate around x=30, rotate around y=0, rotate around z=0, scale=0.5]
    \cube{0}{0}{0} \cube{1}{0}{0} \cube{2+2}{0}{0} \cube{3+2}{0}{0} \cube{4+2}{0}{0}
    
    \cube{0}{1}{0} \cube{1}{1}{0} \cube{2+2}{1}{0} \cube{3+2}{1}{0} 

    \cube{0}{0}{1} \cube{1}{0}{1} \cube{2+2}{0}{1}

    \cube{0}{1}{1} \cube{1}{1}{1} \cube{2+2}{1}{1}

    \cube{0}{0}{2+2} \cube{1}{0}{2+2} \cube{2+2}{0}{2+2}

    \cube{0}{1}{2+2} \cube{1}{1}{2+2}

    \cube{0}{0}{3+2}

\end{tikzpicture}
 \caption{Sections for $\tau$.}
 \label{fig:2}
\end{figure}
For the second section, the only restriction is that there are at most $s^2$ parts. Therefore, it has a generating function of $\frac{1}{(q;q)_{s^2}}$. For the third section, the only restriction is that the length of the parts is at most $s$. Therefore, it has a generating function of $\frac{1}{(q;q)_{s}}$. However, note that the fourth section depends on both the second and third sections. So, we introduce variables to keep track of certain restrictions for the second and third sections.
Note: we set $s>0$ and then add 1 to the final expression for the empty partition. Let $t$ denote the $s^2$th part in the second section. Thus, if the second section has exactly $s^2$ parts, then $t>0$ is its smallest part. Hence, we can remove a rectangle of size $t\times s^2$ from the Ferrers diagram for the second section. The remaining partition has at most $s^2-1$ parts. This has a generating function of $$\frac{q^{ts^2}}{(q;q)_{s^2-1}}.$$ If the section has less than $s^2$ parts, then $t=0$ and the generating function is $$\frac{q^{0s^2}}{(q;q)_{s^2-1}}.$$ Combining the two cases and summing over all values of $t$ yields
$$
\frac{1}{(q;q)_{s^2}}=\sum_{t\geq 0}\frac{q^{ts^2}}{(q;q)_{s^2-1}}.
$$
Note that this also holds by using a geometric series and summing over $t$.
For the third section, we keep track of the number of parts that are equal to $s$. We denote this by $r$. Taking the conjugate, $r$ is the $s$th part. By the same logic as the second section, we have a generating function
$$
\sum_{r\geq 0}\frac{q^{rs}}{(q;q)_{s-1}}=\frac{1}{(q;q)_{s}}.
$$
For the fourth section, the largest part is at most $t$. There can be at most $(s+1)^2-s^2-1$ parts by the definition of the largest cube. However, if the third section has exactly $r$ parts equal to $s$, then the number of parts of the fourth section is at most $r$. Hence, the fourth section has at most $\alpha(r,s)=\min{\{r,2s\}}$ parts. Thus, for the fourth section we have a partition in a box of size $t\times\alpha(r,s)$. The generating function is ${\begin{bmatrix}
        t+\alpha(r,s)\\
        t
    \end{bmatrix}}_{q}, \,\text{where } \alpha(r,s)=\min{\{r,2s\}}.$
Combining these restrictions and summing over all values of $s\geq 1$ and $t,r\geq0$, we obtain a generating function for all unrestricted partitions (except the empty partition). Therefore,
\begin{equation}\label{eq2.2}
    \frac{1}{(q;q)_{\infty}}-1=\sum_{s\geq 1}\sum_{r\geq 0}\sum_{t\geq 0}\frac{q^{s^3+ts^2+rs}}{(q;q)_{s^2-1}(q;q)_{s-1}} {\begin{bmatrix}
        t+\alpha(r,s)\\
        t
    \end{bmatrix}}_{q},
\end{equation}
where $ \alpha(r,s)=\min{\{r,2s\}}$.
Using the fact that 
$$
\sum_{t\geq 0}q^{ts^2} {\begin{bmatrix}
        t+\alpha(r,s)\\
        t
    \end{bmatrix}}_{q}=\frac{1}{(q^{s^2};q)_{\alpha(r,s)+1}},
$$
we obtain by equation (\ref{eq2.2}): 
\begin{align*}
    \frac{1}{(q;q)_{\infty}}-1&=\sum_{s\geq 1}\sum_{r\geq 0}\frac{q^{s^3+rs}}{(q;q)_{s^2-1}(q;q)_{s-1}(q^{s^2};q)_{\alpha(r,s)+1}}\\
    &=\sum_{s \geq 1}\sum_{r\geq0}\frac{q^{s^3+rs}}{(q;q)_{s-1}(q;q)_{s^2+\alpha(r,s)}}, \,\text{where } \alpha(r,s)=\min{\{r,2s\}}.\qedhere
\end{align*}
\end{proof}
Applying the same argument as in Theorem \ref{thmr1} to higher dimensions, we generalise Theorem \ref{thmr1}. 
\begin{theorem}\label{thmr3.3}
Let $k\geq 1$ be an integer. Then, for $|q|<1$, we have
\begin{align}
    \frac{1}{(zq;q)_{\infty}}=&1+\sum_{s \geq 1}\sum_{r\geq0}\frac{z^sq^{s^k+rs}}{(q;q)_{s-1}(zq;q)_{s^{k-1}+\alpha(r,s)}}\label{eq3.5}\\
    =&1+\sum_{s \geq 1}\sum_{r\geq0}\frac{z^{s^{k-1}+r}q^{s^k+rs}}{(zq;q)_{s-1}(q;q)_{s^{k-1}+\alpha(r,s)}},\label{eq3.6}\\
    &\;\;\;\;\text{where } \alpha(r,s)=\min{\{r,(s+1)^{k-1}-s^{k-1}-1\}}.\notag
\end{align}
\end{theorem}
\begin{proof}
Since we already gave a detailed proof utilising the technique for three dimensions, we skip some of the details and explanations.
    Fix $k\geq 3$. Note that $z$ keeps track of the number of parts of a partition. We first prove equation (\ref{eq3.5}). Consider a partition $\lambda=(\lambda_1,\,\lambda_2,\dots,\lambda_j)$ of $n$ with a Durfee modular $k$-cube of side length $s$. Now consider the $s^{k-2}$-modular diagram of $\lambda$. By the definition of the Durfee modular cube, we know that there is an $s\times s$ square in the $s^{k-2}$-modular diagram of $\lambda$. For example, consider a partition $40+35+30+29+28+2$ for $k=4$ that has a Durfee modular 4-cube of side length 3.
    \begin{figure}[H]
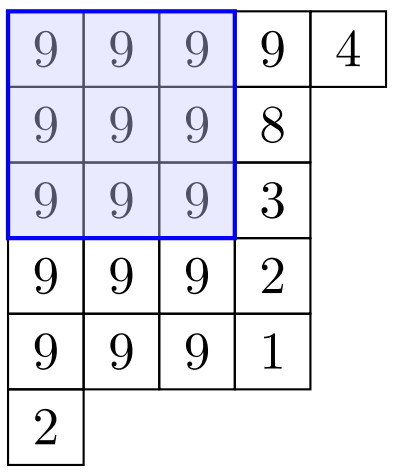

    \centering
        \modgram{9}{40,35,30,29,28,2}{3}
        \caption{The $3^2$-modular diagram for $40+35+30+29+28+2$.}
    \end{figure}
    As in the proof of Theorem \ref{thmr1}, we split the diagram into 4 sections. The $k$-cube has $s$ parts hence contributes $z^sq^{s^k}$. For the partition on the right of the cube, let $r$ denote the $s$th part. The generating function is $$\frac{q^{rs}}{(q;q)_{s-1}}.$$ For the partition below the $k$-cube, let $t$ be the number of parts equal to $s^{k-1}$. We also need to keep track of the number of parts. Thus, the generating function is $$\frac{z^tq^{ts^{k-1}}}{(zq;q)_{s^{k-1}-1}}.$$ Let $\alpha(r,s)=\min{\{r,(s+1)^{k-1}-s^{k-1}-1\}}$. The final remaining partition has a largest part less than or equal to $\alpha(r,s)$ and at most $t$ parts. Applying the same logic as in the proof of Theorem \ref{thmr1}, we obtain
    $$
    \frac{1}{(zq;q)_{\infty}}-1=\sum_{s\geq 1}\sum_{r\geq 0}\sum_{t\geq 0}\frac{z^{s+t}q^{s^k+ts^{k-1}+rs}}{(zq;q)_{s^{k-1}-1}(q;q)_{s-1}} {\begin{bmatrix}
        t+\alpha(r,s)\\
        t
    \end{bmatrix}}_{q}.
    $$
Once again, using the fact that 
$$
\sum_{t\geq 0}(zq^{s^k})^t {\begin{bmatrix}
        t+\alpha(r,s)\\
        t
    \end{bmatrix}}_{q}=\frac{1}{(zq^{s^k};q)_{\alpha(r,s)+1}},
$$
   we obtain the desired result.

   The proof for equation (\ref{eq3.6}) is similar, except that we use the $s^{k-2}$-stacked diagram and the Durfee stacked $k$-cube and then keep track of the number of parts. Break the stacked Ferrers diagram into the same sections as for the proof of equation (\ref{eq3.5}). By the definition of a stacked Durfee cube, we know the cube has $s^{k-1}$ parts equal to $s$ and therefore the generating function for the cube is $z^{s^{k-1}}q^{s^k}$. Let $t$ denote the $s^{k-1}$th part of the partition to the right of the cube. Note that the partition on the right of the cube has at most $s^{k-1}$ parts. Therefore, the generating function is $$
   \frac{q^{ts}}{(q;q)_{s^{k-1}-1}}.
   $$
   Let $r$ denote the number of parts equal to $s$ for the partition below the cube. Thus, the generating function for the partition below the cube is
   $$
   \frac{z^rq^{rs}}{(zq;q)_{s-1}}.
   $$
   The partition in the fourth section has all parts less than or equal to $t$ and fewer than $\alpha(r,s)$ parts. Thus, combining all these generating functions, summing over all possible values of $s,t,r$, and applying the same logic as in the proof of equation (\ref{eq3.5}), the result follows.
\end{proof}
\begin{remark}
    An interesting result to notice is that if we set $k=1$, then Theorem \ref{thmr3.3} simplifies to Theorem \ref{thmp3}. Substitute $z=1$ in equation (\ref{eq3.5}) or (\ref{eq3.6}) to obtain
    $$
    \frac{1}{(q;q)_{\infty}}=1+\sum_{s \geq 1}\sum_{r\geq0}\frac{q^{s^k+rs}}{(q;q)_{s-1}(q;q)_{s^{k-1}+\alpha(r,s)}}.
    $$
    Substitute $k=1$ and $\alpha(r,s)=\min{\{r,(s+1)^{0}-s^{0}-1\}}=\min{\{r,-1\}}=-1$ to see that
    \begin{align*}
        \frac{1}{(q;q)_{\infty}}&=1+\sum_{s \geq 1}\sum_{r\geq0}\frac{q^{s+rs}}{(q;q)_{s-1}(q;q)_{0}}\\
        &=1+\sum_{s \geq 1}\frac{q^{s}}{(q;q)_{s-1}}\sum_{r\geq0}q^{rs}\\
        &=1+\sum_{s \geq 1}\frac{q^{s}}{(q;q)_{s-1}(1-q^s)}\\
        &=\sum_{s \geq 0}\frac{q^{s}}{(q;q)_{s}}.
    \end{align*}
    However, this would imply that in the proof we are using an $s^{-1}$-modular diagram.
    Additionally, notice that setting $k=2$ and $z=1$ in Theorem \ref{thmr3.3} results in the Cauchy identity and setting $k=3$ and $z=1$ results in Theorem \ref{thmr1}. Substitute $z=1$ in equation (\ref{eq3.5}) and notice that for $k=2$, we have $\alpha(r,s)=\min{\{r,(s+1)^{1}-s^{1}-1\}}=\min{\{r,0\}}=0$. Thus,
    \begin{align*}
        \frac{1}{(q;q)_{\infty}}&=1+\sum_{s \geq 1} \sum_{r\geq0}\frac{q^{s^2+rs}}{(q;q)_{s-1}(q;q)_{s^{1}+0}}\\
        &=1+\sum_{s \geq 1}\frac{q^{s^2}}{(q;q)_{s-1}(q;q)_{s}}\sum_{r\geq0}q^{rs}\\
        &=\sum_{s \geq 0}\frac{q^{s^2}}{(q;q)_{s}^2}.
    \end{align*}
\end{remark}
\begin{cor}\label{corr1}
For every integer $k\geq 2$ and $|q|<1$, we have
\begin{align*}
    \frac{1}{(q;q)_{\infty}}=1+&\sum_{s \geq 1}\left(\sum_{r=0}^{\beta(s)}\frac{q^{s^k+rs}}{(q;q)_{s-1}(q;q)_{s^{k-1}+r}}+\frac{q^{s(s+1)^{k-1}}}{(q;q)_{s}(q;q)_{(s+1)^{k-1}-1}}\right),\\
    &\text{ where } \beta(s)=(s+1)^{k-1}-s^{k-1}-1.
\end{align*}
\end{cor}
\begin{proof}
We start by setting $z=1$ in Theorem \ref{thmr3.3}. We split the sum over $r$ into the cases $r\leq \beta(s)$ and $r>\beta(s)$ and apply the geometric series formula to the second case. Thus,
\begin{align*}
    \frac{1}{(q;q)_{\infty}}&=1+\sum_{s \geq 1}\left(\sum_{r=0}^{\beta(s)}\frac{q^{s^k+rs}}{(q;q)_{s-1}(q;q)_{s^{k-1}+r}}+\sum_{r=\beta(s)+1}^{\infty}\frac{q^{s^k+rs}}{(q;q)_{s-1}(q;q)_{s^{k-1}+\beta(s)}}\right)\\
    &=1+\sum_{s \geq 1}\hk*{\sum_{r=0}^{\beta(s)}\frac{q^{s^k+rs}}{(q;q)_{s-1}(q;q)_{s^{k-1}+r}}+\frac{q^{s^k+(\beta(s)+1)s}}{(1-q^s)(q;q)_{s-1}(q;q)_{s^{k-1}+\beta(s)}}}\\
    &=1+\sum_{s \geq 1}\hk*{\sum_{r=0}^{\beta(s)}\frac{q^{s^k+rs}}{(q;q)_{s-1}(q;q)_{s^{k-1}+r}}+\frac{q^{s^k+(\beta(s)+1)s}}{(q;q)_{s}(q;q)_{s^{k-1}+\beta(s)}}}\\
    &=1+\sum_{s \geq 1}\hk*{\sum_{r=0}^{\beta(s)}\frac{q^{s^k+rs}}{(q;q)_{s-1}(q;q)_{s^{k-1}+r}}+\frac{q^{s^k+((s+1)^{k-1}-s^{k-1})s}}{(q;q)_{s}(q;q)_{s^{k-1}+(s+1)^{k-1}-s^{k-1}-1}}}\\
    &=1+\sum_{s \geq 1}\left(\sum_{r=0}^{\beta(s)}\frac{q^{s^k+rs}}{(q;q)_{s-1}(q;q)_{s^{k-1}+r}}+\frac{q^{s(s+1)^{k-1}}}{(q;q)_{s}(q;q)_{(s+1)^{k-1}-1}}\right).\qedhere
\end{align*}
\end{proof}
We now use the same technique to derive another new $q$-identity with higher-powered exponents for partitions into distinct parts.
\begin{theorem}\label{thm3.4}
    For $|q|<1$, we have 
    \begin{equation}
        (-q;q)_\infty =1+\sum_{s\geq1}\frac{q^{s^3+s^2(s^2-1)/2}}{(q;q)_{s^2-1}}\left(\sum_{j=0}^{2s}q^{js^2+js+s^2+j(j+1)/2}\frac{(-q;q)_s}{(q^{s^2};q)_{j+1}}+(-q;q)_{s-1}\right).
    \end{equation}
\end{theorem}
\begin{proof}
    Consider a partition $\lambda$ into distinct parts that has a stacked Durfee cube of side length $s$. Once again, we decompose the $s$-stacked Ferrers diagram for $\lambda$. The Durfee stacked cube has a generating function of $q^{s^3}$. Next, let $\pi_1$ be the partition to the right of the cube with at most $s^2$ distinct parts. 
    We now split the proof into two cases.\\
    \textit{Case 1:} The partition $\pi_1$ has less than $s^2$ parts. Then, $\pi_1$ has exactly $s^2-1$ parts and the generating function is $$\frac{q^{s^2(s^2-1)/2}}{(q;q)_{s^2-1}}.$$ Let $\pi_2$ be the partition in the direction of the $y$-axis in front of the cube. We call the remaining section $\pi_3$. Note that for this case, $\pi_3$ is the empty partition and has a generating function of $1$. Because the last part of the cube is $s$, we know that $\pi_2$ can have a largest part of at most $s-1$. Hence, the generating function for $\pi_2$ is $$(-q;q)_{s-1}.$$
    \textit{Case 2:} The partition $\pi_1$ has exactly $s^2$ parts. Let the smallest part of $\pi_1$ be $t$. Thus, we obtain the generating function
    $$
    \frac{q^{s^2(s^2-1)/2+ts^2}}{(q;q)_{s^2-1}}.
    $$ For this case, we decompose the remaining parts differently. Let $\tau_1$ be the partition consisting of the remaining parts of $\lambda$ below the Durfee cube and with all parts greater than $s$. Let $\tau_2$ be the partition consisting of the remaining parts where each part is at most $s$.
    For example, consider Figure \ref{fig:3.6} below:
    \begin{figure}[H]
        \centering
        \begin{tikzpicture}[scale=0.5]

\draw[thick] (0,0) rectangle (3,3);
\node at (1.5,2) {Durfee};
\node at (1.5,1) {cube};

\draw[thick] (3.4,0) rectangle (5,3);
\node at (4.2,1.5) {$\pi_1$};
\node at (4.2,3.35) {$t$};

\draw[thick]
    (5,0) --
    (7,3) --
    (5,3) --
    cycle;

\draw[thick] (0,-2.4) rectangle (3,-0.4);
\node at (1.5, -1.4) {$\tau_1$};

\draw[thick]
    (0,-2.4-0.3) --
    (3,-2.4-0.3) --
    (0,-4.4-0.3) --
    cycle;
\node at (1.15,-3.25) {$\tau_2$};
\draw[thick]
    (3,-2.4) --
    (3,-0.4) --
    (5,-0.4) --
    cycle;

\end{tikzpicture}
        \caption{Top view of the decomposition.}
        \label{fig:3.6}
    \end{figure}
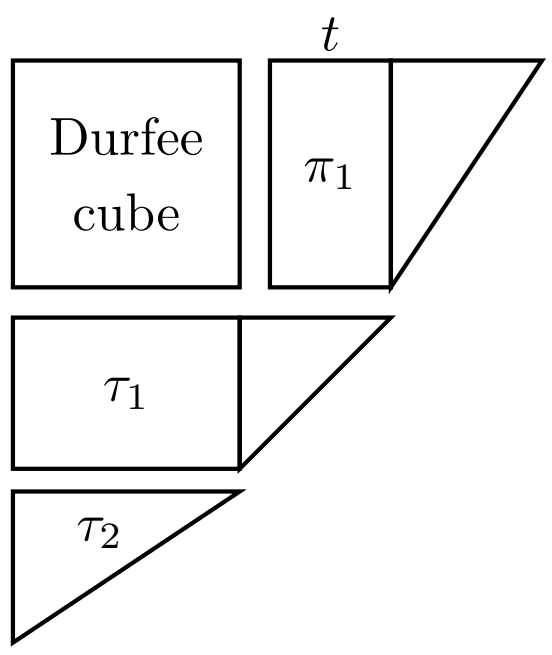
As can be seen in the figure, we can split $\tau_1$ into two partitions: one is a rectangle where the top side is $s$, and the other side is, say, $j$. This contributes the term $q^{sj}$. The other partition is a partition into distinct parts and has its largest part at most $t-1$ (to ensure distinct parts) and needs to have exactly $j$ parts. However, note that $j$ is at most $2s$ to ensure that the Durfee cube has side at most $s$. If we sum over all $j$, we get the generating function
$$
\sum_{j=0}^{2s}q^{js+j(j+1)/2}{\begin{bmatrix}
        t-1\\
        j
    \end{bmatrix}}_{q}.
$$
Finally, for $\tau_2$, since all parts in $\tau_1$ are at least $s+1$, we find that $\tau_2$ is a partition into distinct parts with largest part at most $s$. The generating function of all such partitions is
$$
(-q;q)_s.
$$
Combining all these cases and summing over all possible values of $t$ and $s$, we obtain
\begin{align}\label{eqr3.8}
    (-q;q)_\infty=1+\sum_{s\geq1}\frac{q^{s^3+s^2(s^2-1)/2}}{(q;q)_{s^2-1}}\Biggl(\sum_{t\geq1}\sum_{j=0}^{2s}q^{ts^2+js+j(j+1)/2}{\begin{bmatrix}
        t-1\\
        j
    \end{bmatrix}}_{q}(-q;q)_s\\
    +(-q;q)_{s-1}\Biggr).\notag
\end{align}
Using the $q$-binomial identity $$\sum\limits_{n\geq 0}x^n{\begin{bmatrix}
        n\\
        j
    \end{bmatrix}}_{q}=\frac{x^j}{(x;q)_{j+1}},$$ we can shift the index $t$ to include 0, factor the term $q^{s^2}$ and substitute $x=q^{s^2}$ to find that
    $$
    \sum_{t\geq 1}q^{s^2 t}{\begin{bmatrix} t-1\\ j \end{bmatrix}}_{q}=q^{s^2}\hk*{\frac{q^{js^2}}{(q^{s^2};q)_{j+1}}}.
    $$
    Substituting this into equation (\ref{eqr3.8}) completes the proof.
\end{proof}
We now adjust the technique slightly to derive a new $q$-identity for partitions into distinct odd parts. We accomplish this by using the fact that the number of partitions of $n$ into distinct odd parts is equal to the number of self-conjugate partitions of $n$. Because of the symmetry of self-conjugate partitions, we use the maximal value $s$ such that a partition contains at least $s^2$ parts greater than or equal to $s^2$, instead of Durfee cubes.
\begin{theorem}\label{thmr3.5}
    For $|q|<1$, we have
    \begin{align*}
    (-q;q^2)_\infty1+\sum_{s\geq 1}\sum_{t\geq0}\sum_{r=0}^{\alpha(t,s)}\frac{q^{s^4+2ts^2+r^2}}{(q^2;q^2)_{s^2-1}}{\begin{bmatrix}
        t\\
        r
    \end{bmatrix}}_{q^2},\\
    \text{ where } \alpha(t,s)=\min{(t,2s)}.
\end{align*}
\end{theorem}
\begin{proof}
    Start with a self-conjugate partition $\lambda$ of $n$. Note that the generating function for all such partitions is $(-q;q^2)_\infty$. Let $s$ be the largest integer such that $\lambda$ has at least $s^2$ parts greater than or equal to $s^2$. Consider the Ferrers diagram for $\lambda$. We can draw an $s^2\times s^2$ square inside the Ferrers diagram for $\lambda$. 
    
    Once again, we decompose the diagram into four sections. The first section is the square of side length $s^2$ and contributes the term $q^{s^4}$. The second is the partition to the right of the square, say $\pi_1$, and has at most $s^2$ parts and $s^2$th part, say $t$. The third is the partition below the square, say $\pi_1'$, where all parts are at most $s^2$. Finally, the fourth is the remaining partition, say $\pi_2$, that sits below $\pi_1$ and to the right of $\pi_1'$. 

\begin{figure}[H]
    \centering
\begin{tikzpicture}[
  dot/.style={circle, fill=black, inner sep=1.2pt}, scale=0.7
]

  \def\s{0.58}

  \foreach \r/\len in {1/9, 2/8, 3/7, 4/6, 5/5, 6/4, 7/3, 8/2, 9/1} {
    \foreach \c in {1,...,\len} {
      \node[dot] at ({(\c-1)*\s}, {-(\r-1)*\s}) {};
    }
  }

  \draw[thick, black, rounded corners=2pt]
    ({-0.35*\s}, {0.35*\s}) rectangle ({3.35*\s}, {-3.35*\s});

  \draw[decorate, decoration={brace, amplitude=5pt}, thick, black]
    ({-0.35*\s}, {0.45*\s}) -- ({3.35*\s}, {0.45*\s})
    node[midway, above=7pt, font=\normalsize] {$s^2$};

  \draw[decorate, decoration={brace, amplitude=5pt, mirror}, thick, black]
    ({-0.45*\s}, {0.35*\s}) -- ({-0.45*\s}, {-3.35*\s})
    node[midway, left=7pt, font=\normalsize] {$s^2$};

  \draw[thick, black, dashed, rounded corners=2pt]
    ({3.65*\s}, {0.35*\s}) -- ({8.35*\s}, {0.35*\s}) -- ({8.35*\s}, {-0.35*\s}) --
    ({7.35*\s}, {-0.35*\s}) -- ({7.35*\s}, {-1.35*\s}) --
    ({6.35*\s}, {-1.35*\s}) -- ({6.35*\s}, {-2.35*\s}) --
    ({5.35*\s}, {-2.35*\s}) -- ({5.35*\s}, {-3.35*\s}) -- ({3.65*\s}, {-3.35*\s}) -- cycle;
  \node[black, font=\large] at ({8.9*\s}, {-0.8*\s}) {$\pi_1$};

  \draw[thick, black, dashed, rounded corners=2pt]
    ({-0.35*\s}, {-3.65*\s}) -- ({3.35*\s}, {-3.65*\s}) -- ({3.35*\s}, {-5.35*\s}) --
    ({2.35*\s}, {-5.35*\s}) -- ({2.35*\s}, {-6.35*\s}) --
    ({1.35*\s}, {-6.35*\s}) -- ({1.35*\s}, {-7.35*\s}) --
    ({0.35*\s}, {-7.35*\s}) -- ({0.35*\s}, {-8.35*\s}) -- ({-0.35*\s}, {-8.35*\s}) -- cycle;
  \node[black, font=\large] at ({1.2*\s}, {-8.85*\s}) {$\pi_1'$};

  \node[black, font=\large] at ({4.6*\s}, {-4.6*\s}) {$\pi_2$};

\end{tikzpicture}
\caption{Decomposition of a self-conjugate partition.}
\end{figure}
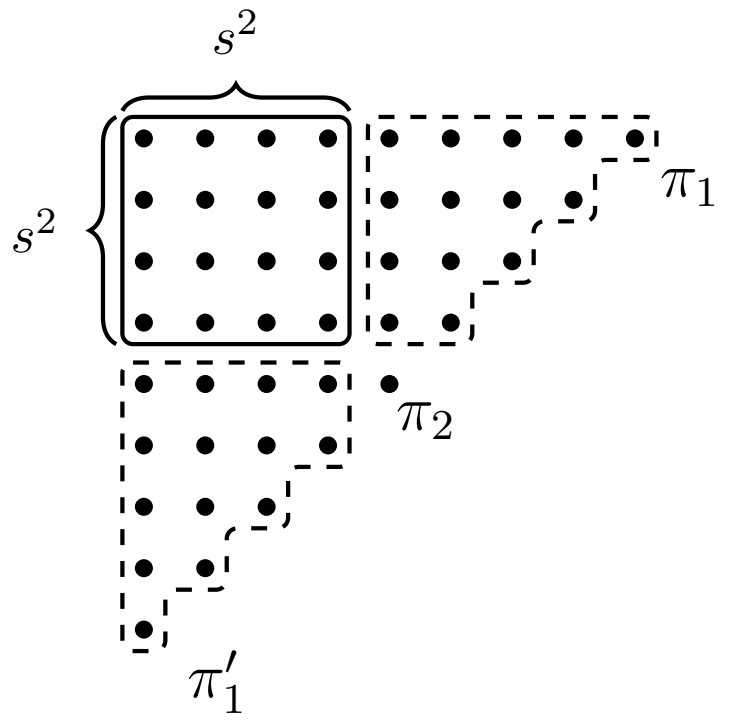
    Since $\lambda$ is a self-conjugate partition, note that $\pi_1$ and $\pi_1'$ are each other's conjugates. Thus, we have that $\pi_1'$ has $t$ parts equal to $s^2$. By adding the conjugate of $\pi_1'$ to $\pi$, we have a partition into even parts with the $s^2$th part equal to $2t$. The generating function for such partitions can be obtained by the same technique as in the above proofs and is
    $$
    \frac{q^{2ts^2}}{(q^2;q^2)_{s^2-1}}.
    $$
    Since $\lambda$ is self-conjugate, we also have that $\pi_2$ is self-conjugate. Additionally, we have the added condition that the largest part of $\pi_2$ is at most $t$ and similarly, the number of parts is at most $t$. Finally, since $s$ is the largest value such that $\lambda$ contains $s^2$ parts greater than or equal to $s^2$, we know that $\pi_2$ does not have $(s+1)^2-s^2-1$ parts greater than or equal to $(s+1)^2-s^2-1$.
    
    Let $r$ be the side of the Durfee square of $\pi_2$. Based on the restrictions, we know that $0\leq r\leq \alpha(t,s)=\min{(t,(s+1)^2-s^2-1)}$. Let $\mu$ be the partition to the right of the Durfee square of $\pi_2$ and $\mu'$ the partition below the Durfee square. Once again, note that $\mu'$ is the conjugate of $\mu$. Note that the parts of $\mu$ can be at most $t-r$ and that there can be at most $r$ parts. If we add the conjugate of $\mu'$ to $\mu$, then we have partitions into even parts that fit in a box of size $2(t-r)\times r$. The generating function for such partitions is
    $$
    {\begin{bmatrix}
        t\\
        r
    \end{bmatrix}}_{q^2}.
    $$
    Finally, combining this with $q^{r^2}$ for the Durfee square of $\pi_2$ and summing over all possible values of $r,s,t$, we obtain the desired result.
\end{proof}
We will now give a generalisation of Theorem \ref{thmc3.4} using Durfee cubes.
\begin{theorem}
    Let $M,N$ be positive integers. Let $k>1$ be an integer. Let $\beta(M,N)=\min{(N,M^{1/(k-1)})}$. Then for $|q|<1$, we have
    \begin{align*}
    \begin{bmatrix}
        M+N\\
        N
    \end{bmatrix}_q=1+\sum_{s=1}^{\beta(N,M)}\sum_{t=0}^{N-s}\sum_{r=0}^{M-s^{k-1}}&q^{s^k+ts^{k-1}+rs}\begin{bmatrix}
        N-t-s-1+s^{k-1}\\
        s^{k-1}-1
    \end{bmatrix}_q\\
    &\times\begin{bmatrix}
        M-r-s^{k-1}-1+s\\
        s-1
    \end{bmatrix}_q\begin{bmatrix}
        t+\alpha(r,s)\\
        t
    \end{bmatrix}_q,
\end{align*}
where $\alpha(r,s)=\min{(r,(s+1)^{k-1}-s^{k-1}-1)}$.
\end{theorem}
\begin{proof}
    Consider a partition $\lambda$ with largest part at most $N$ and at most $M$ parts and a Durfee stacked $k$-cube of side length $s\geq1$. Thus, we have a partition that fits in an $N\times M$ box. Once again, we can decompose the Ferrers diagram for $\lambda$. The generating function for the Durfee $k$-cube is $q^{s^k}$. To the right of the Durfee cube, we have a partition, say $\tau_1$, that fits in a box of size $(N-s)\times s^{k-1}$. Below the Durfee cube, we have a partition, say $\tau_2$, that fits in a box of size $s\times(M-s^{k-1})$. Finally, we have a third partition, say $\tau_3$, that depends on $\tau_1$ and $\tau_2$. Let $t$ be the $s^{k-1}$th part of $\tau_1$, where $t=0$ if $\tau_1$ has less than $s^{k-1}$ parts. Note that $t$ can be in the range $0\leq t\leq N-s$. The left $t$ columns of $\tau_1$ have a generating function of $q^{ts^{k-1}}$. The remaining partition must fit in a box of size $(N-s-t)\times (s^{k-1}-1)$. Thus, the generating function for $\tau_1$ with $s^{k-1}$th part equal to $t$ is
    $$
    q^{ts^{k-1}}\begin{bmatrix}
        N-t-s-1+s^{k-1}\\
        s^{k-1}-1
    \end{bmatrix}_q.
    $$
    For $\tau_2$, we can apply the same logic. Let $r$ denote the number of parts in $\tau_2$ that is equal to $s$. The top $r$ rows of $\tau_2$ have a generating function of $q^{rs}$. The remaining partition must fit into a $(s-1)\times(M-s^{k-1}-r)$ box. The possible range of values for $r$ is $0\leq r\leq M-s^{k-1}$. Therefore, the generating function for $\tau_2$ with $r$ parts equal to $s$ is
    $$
    q^{rs}\begin{bmatrix}
        M-r-s^{k-1}-1+s\\
        s-1
    \end{bmatrix}_q.
    $$
    Finally, $\tau_3$ can have a largest part of at most $t$ and can have at most $r$ parts. However, based on the definition of the Durfee stacked $k$-cube, $\tau_3$ can also have at most $(s+1)^{k-1}-s^{k-1}-1$ parts. Hence, with $\alpha(r,s)$ defined as in the statement of the theorem, we have that $\tau_3$ must fit into a $t\times \alpha(r,s)$ box. Thus, the generating function for $\tau_3$ is 
    $$
    \begin{bmatrix}
        t+\alpha(r,s)\\
        t
    \end{bmatrix}_q.
    $$
    Combining all these generating functions and summing over all possible values of $t$ and $r$, we obtain the generating function for $\lambda$:
    $$
    \sum_{t=0}^{N-s}\sum_{r=0}^{M-s^{k-1}}q^{s^k+ts^{k-1}+rs}\begin{bmatrix}
        N-t-s-1+s^{k-1}\\
        s^{k-1}-1
    \end{bmatrix}_q\begin{bmatrix}
        M-r-s^{k-1}-1+s\\
        s-1
    \end{bmatrix}_q\begin{bmatrix}
        t+\alpha(r,s)\\
        t
    \end{bmatrix}_q.
    $$
    Finally, noting that for $\lambda$ to contain a Durfee stacked $k$-cube of side length $s$ inside and $N\times M$ box, we have that $0\leq s\leq N$ and $o\leq s^{k-1}\leq M$. Thus, with $\beta(M,N)$ defined as in the statement of the theorem, we have $0\leq s\leq \beta(M,N)$. Summing over this range of values for $s$ yields the required result.
\end{proof}
Notice that Theorem \ref{thmr3.5} only worked for even exponents. We now give a theorem for self-conjugate partitions that has a cube in the exponent.
\begin{theorem}
Let $\beta(s)=s^2-s+1$ and $\alpha(s,t)=(\min{(\beta(s)+2t-1,\beta(s+1))}-1)/2$. For $|q|<1$, we have
    $$
    (-q;q^2)_{\infty}=1+\sum_{s\geq1}q^{s^3}\sum_{t\geq0}\frac{q^{2ts}(-q;q^2)_{\alpha(s,t)}}{(q^2;q^2)_{s-1}}.
    $$
\end{theorem}
\begin{proof}
    To prove this theorem, we use the fact that any cube consisting of $s^3$ unit cubes can be reconstructed as a two-dimensional shape consisting of unit squares where there is a square of side length $s$ and a rectangle to the right of the square and another rectangle below the square. This two-dimensional shape is self-conjugate. Therefore, we can use this for self-conjugate partitions. Note that the rectangles next to the square have dimensions $s\times(s^2-s)/2$ and $(s^2-s)/2\times s$.
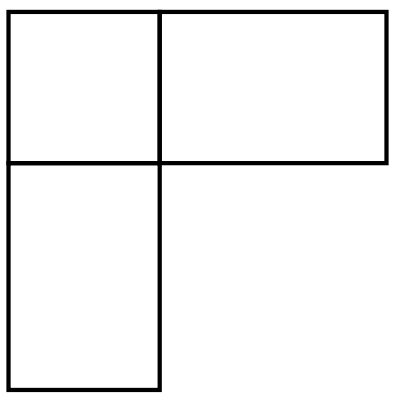
    \begin{figure}[H]
        \centering
\begin{tikzpicture}[scale=0.5]

    \draw[thick] (0,0) rectangle (2,2);

    \draw[thick] (2,0) rectangle (5,2);

    \draw[thick] (0,0) rectangle (2,-3);

\end{tikzpicture}
    \end{figure}
    
   Consider a partition $\lambda$ into distinct odd parts. Let $s$ be the largest value such that $\lambda$ contains at least $s$ parts greater than $\beta(s)=s^2-s+1$. If we place the parts of $\lambda$ to form a self-conjugate Ferrers diagram, then $\lambda$ will contain a shape as the one explained above, which contributes a term of $q^{s^3}$.

    To the right of the right rectangle, we have a partition, say $\pi_1$, into at most $s$ parts with $s$th part, say $t$. Below the lower rectangle, we have a partition, say $\pi_1'$, with largest part at most $s$ and exactly $t$ parts equal to $s$. Note that $\pi_1$ and $\pi_1'$ are conjugates. If we combine the two partitions, we have partitions into even parts with the $s$th part being $2t$. For this, we have a generating function 
    $$\frac{q^{2ts}}{(q^2;q^2)_{s-1}}.$$

    The remaining partition, say $\pi_2$, can be seen as a partition into distinct odd parts with the largest part at most $\beta(s)+2t$. However, since $s$ was maximal, we know that $\lambda$ cannot have $s+1$ parts greater than or equal to $\beta(s+1)$. Thus, we have the added condition that the largest part of $\pi_2$ can be at most $\beta(s+1)-1$. Thus, we set $\alpha(s,t)=(\min{(\beta(s)+2t-1,\beta(s+1))}-1)/2$. Hence, we obtain the generating function
    $$
    (-q;q^2)_{\alpha(s,t)}.
    $$
    Combining all of these generating functions and summing over all possible values of $s$ and $t$, the result follows.
\end{proof}
\section{General Product-Sum Identities Arising from Durfee\\ Rectangle Decompositions}
In the previous section, we proved multiple generalisations of identities using general hypercubes or $k$-cubes. However, this naturally leads to the question of whether these results can be extended to higher-degree polynomials as well. In this section, we will not only extend this technique to certain well-behaved higher-degree polynomials, but to a general class of functions that satisfy certain conditions.

We start by defining a set of functions $\mathcal{P}$. Fix a value $c\in\mathbb{Z}\cup\brac*{-\infty}$ and let $A=\brac*{z\in\mathbb{Z}\,|\,z\geq c}$. Then
\begin{align*}
    \mathcal{P}:=\{f:A\to\mathbb{Z}\,|\,&\text{there exists a constant }c_1 \text{ such that } f(x)>0 \text{ for all } x\geq c_1 \text{ and }\\
    &\text{there exists a constant }c_2 \text{ such that } f(x) \text{ is non-decreasing}\\
    &\text{for all } x\geq c_2\}.
\end{align*}
Some examples of such functions:
\begin{itemize}
    \item polynomials where the coefficient of the term with the highest power is positive, and all coefficients are integers, e.g. $f(x)=3x^7+5x^5-11x^3-x+17$,
    \item functions that generate polygonal numbers, e.g., $f(x)=x(x+1)/2$,
    \item $f(x)=x!$,
    \item exponential functions such as $f(x)=2^x$,
    \item the partition counting function $p(n)$ or the prime counting function $\pi(x)$.
    \item oscillating functions such as $f(x)=x^2+\sin(\pi x)$,
    \item functions with floor or ceiling, e.g., $f(x)=\floor*{e^x}$.
\end{itemize}
Note that in the following theorem, we consider an empty sum to be equal to 1.
\begin{theorem}\label{thmr3.7}
Let $f_1,f_2\in\mathcal{P}$ be functions such that $f_1(x)f_2(x)$ is strictly increasing for all $x$ greater than or equal to some constant $c$. 
We define the following functions:
\begin{itemize}
    \item $\beta_1(s)=f_1(s+1)-f_1(s)-1,\,\beta_2(s)=f_2(s+1)-f_2(s)-1$,
    \item $\gamma(s)=|\beta_2(s)-\beta_1(s)|,\,\kappa(r,s)=\max{(0,r-\gamma(s))}$,
    \item $\alpha(r,s,t)=\min{(\beta_1(s),t,\kappa(r,s))},\, \alpha_2(r,s)=\min{(r,\gamma(s))}$,
    \item $\tau_1=|f_1(a)-f_2(a)|,\, \tau_2=\min{(f_1(a)-1,f_2(a)-1)}$.
\end{itemize}
Based on the definition $f$, we know that there exists an $a\in\mathbb{Z}$ such that $f_1(x),\,f_2(x)$ is strictly positive and non decreasing and $f_1(x)f_2(x)$ is strictly increasing for all $x\geq a$. Then, assuming $\beta_2(s)\geq \beta_1(s)$ for all $s\geq a$, we have
    \begin{align*}
    \frac{1}{(zq;q)_{\infty}}=&\sum_{s\geq a}\frac{z^{f_2(s)}q^{f_1(s)f_2(s)}}{(zq;q)_{f_1(s)-1}(q;q)_{f_2(s)-1}}\sum_{t\geq0}q^{tf_2(s)}\sum_{r\geq0}z^{r}q^{rf_1(s)}\\
    &\times\sum_{i=0}^{\alpha(r,s,t)}q^{i(i+\gamma(s))}{\begin{bmatrix}
        t+\alpha_2(r,s)\\
        i+\alpha_2(r,s)
    \end{bmatrix}}_{q}{\begin{bmatrix}
        \kappa(r,s)\\
        i
    \end{bmatrix}}_{q}\\
    &+\sum_{s=0}^{\tau_2}\frac{z^s q^{s(s+\tau_1)}}{(q;q)_{s}(zq;q)_{s+\tau_1}}.
\end{align*}
\end{theorem}
\begin{proof}
    Let $f_1,f_2\in \mathcal{P}$ satisfy all the conditions of the theorem. Define $\beta_1,\beta_2,\dots,\tau_1,\tau_2$ as in the theorem. We have that $f_1(s)$ and $f_2(s)$ are strictly positive for all $s\geq a$.  

    Consider a partition $\lambda$ of $n$. We are faced with two possibilities: \textit{Case 1:} there exists an $s$ such that $\lambda$ contains a rectangle of dimensions $f_1(s)\times f_2(s)$ and \textit{Case 2:} $\lambda$ does not contain a rectangle for any values of $s\geq a$.\\
    \textit{Case 1:} since $f_1f_2$ is strictly increasing, we know that there exists a largest (hence unique) value $s$ such that $\lambda$ contains a rectangle with dimensions $f_1(s)\times f_2(s)$. Note that $f_1$ will be the side length of the rectangle reading from left to right and $f_2$ will be the side length from top to bottom. We decompose $\lambda$ into four sections again. The first section is the rectangle and contributes the term $$z^{f_2(s)}q^{f_1(s)f_2(s)}$$ where $z$ tracks the number of parts. Next, we have a partition to the right of the rectangle, say $\pi_1$ with at most $f_2(s)$ parts, and a partition below the rectangle, say $\pi_2$ with largest part at most $f_1(s)$. We denote the partition consisting of all remaining parts by $\pi_3$. Let $t$ be the number $f_2(s)$ part of $\pi_1$ and $r$ be the number of parts equal to $f_1(s)$ of $\pi_2$. Once again, the generating function for $\pi_1$ is $$
    \frac{q^{tf_2(s)}}{(q;q)_{f_2(s)-1}}.
    $$ For $\pi_2$, we need to keep track of the number of parts and find a generating function of
    $$
    \frac{z^rq^{rf_1(s)}}{(zq;q)_{f_1(s)-1}}.
    $$
\begin{figure}[H]
    \centering
    \begin{tikzpicture}[scale=0.5]


\draw (0,0) rectangle (4,4);

\draw (4,0) rectangle (6.5,4);

\draw (4,-2) rectangle (6,-0);

\draw (0,-4) rectangle (4,0);

\draw[<->] (0,4.5) -- (4,4.5);
\node at (2,4.9) {$f_1(s)$};

\draw[<->] (4,4.5) -- (6.5,4.5);
\node at (5.25,4.9) {$t$};

\draw[<->] (-0.8,0) -- (-0.8,4);
\node[left] at (-0.8,2) {$f_2(s)$};

\draw[<->] (-0.8,-4) -- (-0.8,0);
\node[left] at (-0.8,-2) {$r$};

\draw[<->] (4.1,-0.2-2.6) -- (6.0,-0.2-2.6);
\node[above] at (5.05,-0.2-2.65) {$i$};

\draw[<->] (6.5,-2) -- (6.5,-0.09);
\node[right] at (6.5,-1) {$i+\gamma(s)$};

\end{tikzpicture}
    \caption{Skeleton of decomposition of $\lambda$.}
\end{figure}
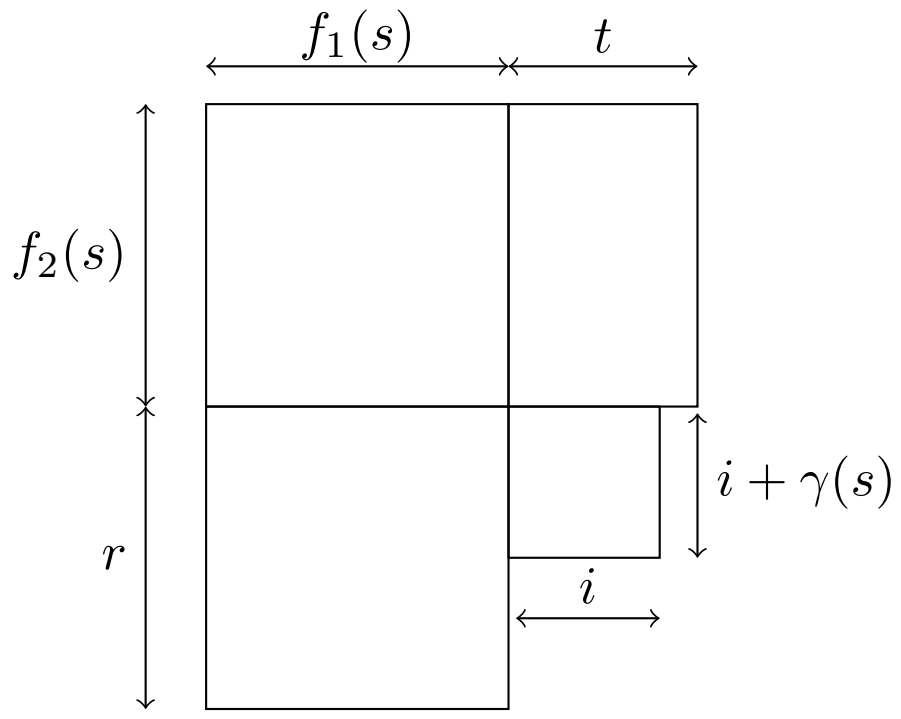
    It remains to find a generating function for $\pi_3$. We know that $\pi_3$ can have a largest part of at most $t$ and can have at most $r$ parts. Additionally, since $s$ is the maximal value, we know that $\lambda$ cannot contain a rectangle with dimensions $f_1(s+1)\times f_2(s+1)$. Therefore, $\pi_3$ can contain at most $\beta_2(s)$ parts greater than or equal to $\beta_1(s)$. Thus, we need to find the generating function for partitions that have a Durfee rectangle of dimensions $i\times (i+\gamma(s))$ for values of $i\leq \beta_1(s)$ (since $\beta_1(s)\leq \beta_2(s)$ by assumption) and fit in a box of dimensions $t\times r$. However, since the largest part of $\pi_3$ is at most $t$, we have that $i\leq t$. Additionally, since we have at most $r$ parts, we have that $i+\gamma(s)\leq r$, thus $i\leq r-\gamma(s)$. Therefore, we have that $$0\leq i\leq \alpha(r,s,t).$$ 
    
    Consider $\pi_3$ that has a Durfee rectangle $i\times (i+\gamma(s))$ for a fixed $i$. The Durfee rectangle contributes the generating function $$q^{i(i+\gamma(s))}.$$ To the right of the Durfee rectangle, we have a partition, say $\mu_1$, and below the Durfee rectangle, we have a partition, say $\mu_2$. For $\mu_1$, we have a partition where the largest part is at most $t-i$. If $\gamma(s)\leq r$, then the number of parts is at most $i+\gamma(s)$. However, if $\gamma(s)>r$, then $\pi_3$ does not contain a Durfee rectangle, hence $i=0$ and the number of parts of $\mu_1=\pi_3$ is at most $r=i+r$. Thus, the number of parts of $\mu_1$ is at most $i+\alpha_2(r,s)$. Hence, we find the generating function $$
    {\begin{bmatrix}
        t+\alpha_2(r,s)\\
        i+\alpha_2(r,s)
    \end{bmatrix}}_{q}.
    $$ For $\mu_2$, we have a partition where the largest part is at most $i$. If $\gamma(s)\leq r$, then the number of parts is at most $r-(i+\gamma(s))$. If $\gamma(s)>r$, then again, we have that $i=0$, thus we have the empty partition with generating function 1. Hence, we obtain the generating function $${\begin{bmatrix}
        \kappa(r,s)\\
        i
    \end{bmatrix}}_{q}.$$

    Finally, combining all the generating functions and summing over all possible values of $s,t,r$, and $i$, we obtain the first quadruple sum on the right-hand side of the theorem.\\
   \textit{ Case 2:} Consider a partition $\lambda$ that does not contain a rectangle with dimensions $f_1(a)\times f_2(a)$. Therefore, we know that $\lambda$ contains a Durfee rectangle with dimensions $s\times (s+\tau_1)$ where $0\leq s\leq\tau_2$. Obtaining the generating function for this is straightforward, and summing over all $s$ gives the second summation in the theorem.

   Combining \textit{Case 1} and \textit{Case 2} completes the proof.
\end{proof}

\begin{cor}
    For a fixed positive integer $k$, and $|q|<1$, we have
    $$\frac{1}{(q;q)_{\infty}}=1+\sum_{s\geq1}\sum_{r\geq0}\sum_{t\geq0}\frac{z^{s^k+r}q^{s^{2k}+(r+t)s^k}}{(zq;q)_{s^k-1}(q;q)_{s^k-1}}\sum_{i=0}^{\alpha(r,s,t)}q^{i^2}{\begin{bmatrix}
        t\\
        i
    \end{bmatrix}}_{q}{\begin{bmatrix}
        r\\
        i
    \end{bmatrix}}_{q},$$
    where $\alpha(r,s,t)=\min((s+1)^k-s^k-1,t,r)$.
\end{cor}
\begin{proof}
    We set $f_1(s)=s^k=f_2(s)$. Note that this satisfies all the required conditions of Theorem \ref{thmr3.7} with $a=1$. Consider
    \begin{itemize}
        \item $\beta_1(s)=(s+1)^k-s^k-1=\beta_2(s)$,
        \item $\gamma(s)=0$, $\kappa(r,s)=r$,
        \item $\alpha(r,s,t)=\min(\beta_1(s),t,r)$, $\alpha_2(r,s)=0$,
        \item $\tau_1=\tau_2=0$.
    \end{itemize}
    Therefore, by Theorem \ref{thmr3.7}, the result follows.
\end{proof}
In the following theorem, we consider the $q$-binomial to be 0 whenever one of the constraints is negative. Therefore, we will set the generating function for the number of partitions that can fit in a box where at least one of the side lengths of the box is negative to be 0.
\begin{theorem}
Let $k$ be a positive integer. Let $\alpha(s)=(s+1)^k-s^{k}-1$. Then the generating function for the number of partitions of $n$ into distinct parts where $z$ tracks the number of parts is given by
    \begin{align*}
    (-zq;q)_\infty=1&+\sum_{s\geq 1}\frac{z^{s^{k}}q^{s^{2k}+s^k(s^k-1)/2}(-zq;q)_{s^k-1}}{(q;q)_{s^k-1}}\\
    &+\sum_{s\geq 1}\sum_{t\geq1}\frac{z^{s^k}q^{s^{2k}+s^k(s^k-1)/2+ts^k}}{(q;q)_{s^k-1}}\\
    &\;\;\;\times \sum_{j=0}^{\alpha(s)}z^{j}\begin{aligned}[t]
    &\left(q^{j(j-1)/2+j^2+js^k}{\begin{bmatrix}
        t-1-j\\
        j-1
    \end{bmatrix}}_{q}(-zq;q)_{s^k+j-1}\right.\\
    &\,\,+\Biggl. q^{j(j+1)/2+j^2+js^k}{\begin{bmatrix}
        t-1-j\\
        j
    \end{bmatrix}}_{q}(-zq;q)_{s^k+j}\Biggr).
    \end{aligned}\\
\end{align*}
\end{theorem}
\begin{proof}
     Fix $k>0$. Let $\lambda$ be a partition into distinct parts. Let $s$ be the largest integer such that the Ferrers diagram of $\lambda$ contains an $s^k\times s^k$ square. Recall the proof of Theorem \ref{thm3.4}. Let $\pi_1$ be the partition to the right of the $s^k\times s^k$ square consisting of at most $s^k$ parts. Let $t$ be the $s^k$th part of $\pi_1$. We have two cases.\\
     \textit{Case 1:} if $t=0$. Then, $\pi_1$ has exactly $s^k-1$ parts and has a generating function $$
     \frac{q^{s^k(s^k-1)/2}}{(q;q)_{s^k-1}}.
     $$
     The square has exactly $s^k$ parts. Therefore, we need to track the number of parts. Thus, we obtain a generating function $$
     z^{s^k}q^{s^{2k}}.
     $$ Finally, below the square, we have a partition into distinct parts with largest part at most $s^k-1$. Additionally, we need to track the number of parts. Thus, the generating function is $(-zq;q)_{s^k-1}$. Combining the three generating functions yields
     $$
     \frac{z^{s^{k}}q^{s^{2k}+s^k(s^k-1)/2}(-zq;q)_{s^k-1}}{(q;q)_{s^k-1}}.
     $$\\
     \textit{Case 2:} if $t>0$. Using the same techniques as the previous proofs, the generating function for the square and $\pi_1$ is $$\frac{z^{s^k}q^{s^{2k}+s^k(s^k-1)/2+ts^k}}{(q;q)_{s^k-1}}.$$

     Consider the partition below the $s^k\times s^k$ square and $\pi_1$, say $\pi_2$. Note that $\pi_2$ is a partition into distinct parts with largest part at most $s^k+t-1$. However, there are additional restrictions on $\pi_2$ to ensure the maximality of $s$. Therefore, we will further decompose $\pi_2$. Let $\tau_1$ be the partition where all parts of $\tau_1$ are strictly below $\pi_1$. Let $j$ denote the side length of the largest square that fits in the partition $\tau_1$. Because of the maximality of $s$, we have $j\leq\alpha(s)$. To the left of the $j\times j$ square, we have an $s^k\times j$ rectangle (with j parts). Therefore, the square and rectangle give us the terms $z^jq^{j^2+js^k}$. We now have two cases again.\\
     \textit{Case A:} the $j$th part of $\tau_1$ is $j$. Then, we have a triangle with $j-1$ parts to the right of the $j\times j$ square which contributes $q^{j(j-1)/2}$. To the right of this triangle, we have a partition that can fit in a box of dimensions $(t-1-j-(j-1))\times(j-1)$. Therefore, we obtain a generating function of $$
     {\begin{bmatrix}
        t-1-j\\
        j-1
    \end{bmatrix}}_{q}.
     $$ Note that if one of the dimensions of the box is negative, then we set the generating function equal to 0. Finally, below the $s^k\times j$ rectangle and $j\times j$ square, we have a partition into distinct parts with largest part at most $s^k+j-1$ where $z$ tracks the number of parts. This has a generating function of $(-zq;q)_{s^k+j-1}$. Therefore, the generating function for \textit{Case A} is $$
     z^{j}q^{j(j-1)/2+j^2+js^k}{\begin{bmatrix}
        t-1-j\\
        j-1
    \end{bmatrix}}_{q}(-zq;q)_{s^k+j-1}.
     $$\\
     \textit{Case B:} if the $j$th part of $\tau_1$ is greater than $j$. Then, we have a triangle with $j$ parts to the right of the $j\times j$ square, which contributes the term $q^{j(j+1)/2}$. Therefore, the partition to the right of the triangle must fit in a box of dimensions $(t-1-j-j)\times j$. Thus, the generating function is $$ {\begin{bmatrix}
        t-1-j\\
        j
    \end{bmatrix}}_{q}.$$ Finally, below the square and rectangle, we have a partition into distinct parts with the largest part at most $s^k+j$. Therefore, the generating function for \textit{Case B} is $$
    z^jq^{j(j+1)/2+j^2+js^k}{\begin{bmatrix}
        t-1-j\\
        j
    \end{bmatrix}}_{q}(-zq;q)_{s^k+j}.
    $$
    Combining all possible cases and summing over all possible values of $s,\,t$ and $j$, the result follows.
\end{proof}
Setting $z=-1$ in the above theorem, we get an interesting corollary.
\begin{cor}
    Let $k$ be a positive integer. Let $\alpha(s)=(s+1)^k-s^{k}-1$. Then
    \begin{align*}
    (q;q)_\infty=1&+\sum_{s\geq 1}(-1)^{s^{k}}q^{s^k(3s^k-1)/2}\\
    &+\sum_{s\geq 1}\sum_{t\geq1}\frac{(-1)^{s^k}q^{s^k(3s^k-1)/2+ts^k}}{(q;q)_{s^k-1}}\\
    &\;\;\;\times \sum_{j=0}^{\alpha(s)}(-1)^{j}\begin{aligned}[t]
    &\left(q^{j(j-1)/2+j^2+js^k}{\begin{bmatrix}
        t-1-j\\
        j-1
    \end{bmatrix}}_{q}(q;q)_{s^k+j-1}\right.\\
    &\,\,+\Biggl. q^{j(j+1)/2+j^2+js^k}{\begin{bmatrix}
        t-1-j\\
        j
    \end{bmatrix}}_{q}(q;q)_{s^k+j}\Biggr).
    \end{aligned}\\
\end{align*}
\end{cor}
Note that by setting $k=1$ in the above corollary, we can derive Euler's pentagonal number theorem.
\begin{theorem}
    Let $f\in \mathcal{P}$ be a function such that $f(x)$ is strictly positive and strictly increasing for all $x\geq a$. Let $\alpha(t,s)=\min{(t,f(s+1)-f(s)-1)}$. For $|q|<1$, we have
    \begin{align*}
    (-q;q^2)_\infty=\sum_{s=0}^{f(a)-1}\frac{q^{s^2}}{(q^2;q^2)_s}+\sum_{s\geq a}q^{f(s)^2}\sum_{t\geq0}\frac{q^{2tf(s)}}{(q^2;q^2)_{f(s)-1}}\sum_{r=0}^{\alpha(t,s)}q^{r^2}{\begin{bmatrix}
        t\\
        r
    \end{bmatrix}}_{q^2}.
\end{align*}
\end{theorem}
\begin{proof}
    Consider a partition $\lambda$ into distinct odd parts. If $\lambda$ contains a square of side length $f(s)$, then by applying the same argument as in the proof of Theorem \ref{thmr3.5}, we find the right-hand triple sum. If $\lambda$ does not contain a square of side length $f(s)$ for some value $s\geq a$, then we have a partition into distinct odd parts with a Durfee square that has a side length of at most $f(a)-1$. Obtaining the generating function for this is straightforward. Combining the two cases, the result follows.
\end{proof}
\bibliography{bib}
\end{document}